\documentclass[11pt,a4paper]{amsart}
\usepackage{amsmath}
\usepackage{amssymb}
\usepackage{amsthm}
\usepackage{amsfonts}
\usepackage{microtype}
\usepackage{mathrsfs}
\usepackage{graphicx}
\usepackage{color}
\usepackage{latexsym}
\usepackage{rotating}
\usepackage{xspace}
\usepackage[all]{xy}
\usepackage{longtable}
\usepackage{leftidx}
\usepackage{mathtools}
\usepackage{titlesec}
\usepackage{tikz}
\usetikzlibrary{calc}
\usetikzlibrary{arrows}
\usetikzlibrary{decorations.markings}
\usepackage{geometry}
\usepackage{amscd}
\usepackage{latexsym}
\usepackage{pstcol,pst-plot,pst-3d}
\usepackage{multicol}
\usepackage{amsmath}
\usepackage{amssymb}
\usepackage{amsthm}
\usepackage{pgfplots}
\psset{unit=0.7cm,linewidth=0.8pt,arrowsize=2.5pt 4}

\newpsstyle{fatline}{linewidth=1.5pt}
\newpsstyle{fyp}{fillstyle=solid,fillcolor=verylight}
\definecolor{verylight}{gray}{0.97}
\definecolor{light}{gray}{0.9}
\definecolor{medium}{gray}{0.85}

\newtheorem{theorem}{Theorem}[section]
\numberwithin{equation}{theorem}
\newtheorem{lemma}[theorem]{Lemma}
\newtheorem{proposition}[theorem]{Proposition}
\newtheorem{corollary}[theorem]{Corollary}

\theoremstyle{definition}
\newtheorem{definition}[theorem]{Definition}
\theoremstyle{fact}
\newtheorem{fact}[theorem]{Fact}
\newtheorem{example}[theorem]{Example}

\newtheorem{remark}[theorem]{Remark}

\theoremstyle{conjecture}

\newcommand{\Ass}{\operatorname{Ass}}
\newcommand{\im}{\operatorname{im}}

\newcommand{\grade}{\operatorname{grade}}

\newcommand{\cd}{\operatorname{cd}}

\newcommand{\f}{\operatorname{f}}

\newcommand{\Ext}{\operatorname{Ext}}

\newcommand{\Tor}{\operatorname{Tor}}

\newcommand{\depth}{\operatorname{depth}}

\newcommand{\fm}{\frak{m}}
\newcommand{\fp}{\frak{p}}
\newcommand{\fq}{\frak{q}}
\newcommand{\fa}{\frak{a}}
\newcommand{\fb}{\frak{b}}

\newcommand{\fn}{\frak{n}}

\newcommand{\suchthat}{\;\ifnum\currentgrouptype=16 \middle\fi|\;}

\newenvironment{prf}[1][Proof]{\begin{proof}[\bf #1]}{\end{proof}}

\makeatletter
\newcommand{\holim@}[2]{%
\vtop{\m@th\ialign{##\cr
\hfil$#1\operator@font holim$\hfil\cr
\noalign{\nointerlineskip\kern1.5\ex@}#2\cr
\noalign{\nointerlineskip\kern-\ex@}\cr}}%
}
\newcommand{\holim}{%
\mathop{\mathpalette\holim@{\rightarrowfill@\textstyle}}\nmlimits@
}
\makeatother

\makeatletter
\def\@secnumfont{\bfseries}
\def\section{\@startsection{section}{1}%
\z@{.7\linespacing\@plus\linespacing}{.5\linespacing}%
{\normalfont\Large\bfseries\filcenter}}
\def\subsection{\@startsection{subsection}{2}%
\z@{.5\linespacing\@plus.7\linespacing}{-.5em}%
{\normalfont\large\bfseries}}
\makeatother
\begin{document}

\author[A. Rahimi]
{Ahad Rahimi}

\title[Algebraic properties of tensor products ...]
{Algebraic properties of tensor product of modules over a field}

\address{A. Rahimi, Department of Mathematics, Razi University, Kermanshah, Iran.}
\email{ahad.rahimi@razi.ac.ir}
\subjclass[2020]{13C14; 13C15; 18G40; 13D45.}

\keywords { Cohen--Macaulay module; cohomological dimension;  dimension filtration; local cohomology;
   finiteness dimension;  sequentially Cohen--Macaulay module; spectral sequence.}
\maketitle
\begin{center}
\textit{Dedicated to the memory of Professor J\"urgen Herzog,\\ 
whose contributions to the field of mathematics\\
  will always be remembered.}
\end{center}

\begin{abstract}
Let $A$ and $B$ be commutative Noetherian algebras over an arbitrary field $\Bbbk$ such that $A \otimes_\Bbbk B$ is Noetherian. We consider ideals $I$ and $J$ of $A$ and $B$, respectively, as well as nonzero finitely generated modules $L$ and $N$ over $A$ and $B$, respectively. In this paper, we investigate certain algebraic properties of the $A \otimes_\Bbbk B$-module $L\otimes_{\Bbbk} N$, which are often inherited from the properties of the $A$-module $L$ and the $B$-module $N$. Specifically, we provide characterizations for the Cohen-Macaulayness, generalized Cohen-Macaulayness, and sequentially Cohen-Macaulayness of $L\otimes_{\Bbbk} N$ with respect to the ideal $I \otimes_\Bbbk B + A \otimes_\Bbbk J$, in terms of the corresponding properties for $L$ and $N$ with respect to $I$ and $J$, respectively.
\end{abstract}

\maketitle

\tableofcontents

\section{Introduction}
Let $A$ and $B$ be commutative Noetherian algebras over a field $\Bbbk$, with $A \otimes_\Bbbk B$ also being Noetherian. Let $I$ and $J$ be ideals of $A$ and $B$, respectively. Then,   $I\otimes_\Bbbk B$ and $A\otimes_\Bbbk J$ are also ideals of $A\otimes_{\Bbbk} B$. We define $\mathcal{I}= I\otimes_\Bbbk B + A \otimes_\Bbbk J$. 
We can also consider $A \otimes_\Bbbk B$ as a local ring with the unique maximal ideal $\mathfrak{M} = \fm \otimes_\Bbbk B + A \otimes_\Bbbk \fn$, where $(A, \fm)$ and $(B, \fn)$ are local algebras over the same field $\Bbbk$. In this case, $A$ and $B$ share $\Bbbk$ as the common residue field, and either $A$ or $B$ is algebraic over $\Bbbk$. This fact is mentioned in Remark \ref{2.8}.
 Moreover, let $L$ and $N$ be nonzero finitely generated modules over $A$ and $B$, respectively. In this paper, we investigate the algebraic properties of the $A\otimes_{\Bbbk} B$-module $L\otimes_{\Bbbk} N$, which are connected to those of $L$ and $N.$  In fact, we provide characterizations for the Cohen-Macaulayness, generalized Cohen-Macaulayness, and sequentially Cohen-Macaulayness of $L\otimes_{\Bbbk} N$ with respect to $\mathcal{I}$, in terms of the corresponding properties for $L$ and $N$ with respect to $I$ and $J$, respectively. It is important to note that in order to study these algebraic properties for $L\otimes_{\Bbbk} N$, we require the assumption that $A\otimes_{\Bbbk} B$ is Noetherian. In general, if $A$ and $B$ are Noetherian algebras, then $A\otimes_{\Bbbk} B$ is not necessarily Noetherian, even when $A$ and $B$ are field extensions of $\Bbbk$, as shown in \cite{V}.
 
Let us begin by discussing Cohen-Macaulayness.
Let $R$ be a Noetherian ring, $\fa$ a proper ideal of $R$, and $M$ a finitely generated $R$-module. We say that $M$ is Cohen-Macaulay with respect to $\fa$ if either $M=\fa M$ or $M\neq \fa M$ and $\grade(\fa, M)=\cd(\fa, M)$ (see \cite{R2}). Here, $\cd(\fa, M)$ denotes the cohomological dimension of $M$ with respect to $\fa$, which is the supremum of all integers  $i$ for which $\text{H}^i_\fa(M)\neq 0$.
We observe that if $(R, \fm')$ is a local ring, then $M$ is Cohen-Macaulay if and only if $M$ is Cohen-Macaulay with respect to $\fm'$.

The purpose of Section 2 is to explain how $\grade(\mathcal{I}, L\otimes_{\Bbbk} N)$ and $\cd(\mathcal{I}, L\otimes_{\Bbbk} N)$ can be expressed in terms of the
 corresponding quantities of $L$ and $N$ with respect to $I$ and $J$. Indeed, in Theorem \ref{2.6}, we demonstrate that
  $\grade(\mathcal{I}, L\otimes_{\Bbbk} N)=\grade(I, L)+\grade(J, N)$ and $\cd(\mathcal{I}, L\otimes_{\Bbbk} N)=\cd(I, L)+\cd(J, N).$
   Thus, we establish that $L\otimes_{\Bbbk} N$ is Cohen--Macaulay with respect to $\mathcal{I}$ if and only if $L$
    and $N$ are Cohen--Macaulay with respect to $I$ and $J$, respectively. Consequently, we obtain \cite[Theorem 1.1]{BK}(b) and \cite[Theorem 2.1]{BK} as corollaries.

Let $R$ be a Noetherian ring, $\fa$ an ideal of $R$, and $M$ a finitely generated $R$-module. We say $M$ is \textit{generalized Cohen--Macaulay with respect to $\fa$} if $\cd(\fa, M)\leq 0$ or $\f_\fa(M)=\cd(\fa, M)$   (see \cite{DGTZ}). Here, $\f_\fa(M)$ denotes the \textit{finite dimension of $M$ with respect to $\fa$}, which is defined as the infimum of non-negative integers $i$ for which $\text{H}_\fa^i(M)$ is not finitely generated. Note that if $(R, \fm')$ is a local ring, then $M$ is generalized Cohen--Macaulay if and only if $M$ is generalized Cohen--Macaulay with respect to $\fm'$.

In Section 3, we provide a description of $\f_{\mathcal{I}}(L\otimes_{\Bbbk} N)$ in terms of the finite dimension and grade of $L$ and $N$ with respect to $I$ and $J$, respectively. Theorem \ref{3.2} states that
\[
\f_{\mathcal{I}}(L\otimes_{\Bbbk} N)= \min\{\grade(I, L)+\f_J(N), \f_I(L)+\grade(J, N)\}.
\]
In particular, in the case of an ordinary local ring, we have
\[
\f_{\mathfrak{M}}(L\otimes_{\Bbbk} N)= \min\{\depth_A L+\f_\fn(N), \f_\fm (L)+\depth_B N \}.
\]
Let  $\cd(I, L)\geq 1$ and $\cd(J, N)\geq 1$. Using this result, we can demonstrate that $L\otimes_{\Bbbk} N$ is generalized Cohen-Macaulay with respect to $\mathcal{I}$ if and only if it is Cohen-Macaulay with respect to $\mathcal{I}$. This, in turn, is equivalent to $L$ and $N$ being Cohen--Macaulay with respect to $I$ and $J$, respectively. Moreover, if one of these equivalent conditions holds, then we have $\f_{\mathcal{I}}(L\otimes_{\Bbbk} N)=\f_I(L)+\f_J(N)$, see Proposition \ref{3.4}.
In particular, for the local case:  Let  $\dim_A L$ and $\dim_B N$ be positive. Then, $L\otimes_{\Bbbk} N$ is generalized Cohen--Macaulay if and only if it is Cohen--Macaulay, which is further equivalent to $L$ and $N$ both being Cohen--Macaulay. Moreover, we have $\f_{\mathfrak{M}}(L\otimes_{\Bbbk} N)=\f_\fm(L)+\f_\fn(N)$, see a graded version of this result in \cite[Theorem 2.6]{STY}. Example \ref{3.7} provides clarity on this result.

Let $R$ be a Noetherian ring, $\fa$ an ideal of $R$, and $M$ a finitely generated $R$-module. We say that $M$ has a \textit{Cohen-Macaulay filtration with respect to $\fa$} if there exists a finite filtration $\mathcal{F}$: $0=M_0\varsubsetneq M_1 \varsubsetneq \dots \varsubsetneq M_r=M$ of $M$ by submodules $M_i$, such that each quotient $M_i/M_{i-1}$ is Cohen-Macaulay with respect to $\fa$, and $0 \leq \cd(\fa, M_1/M_0)<\cd(\fa, M_2/M_1)< \dots< \cd(\fa, M_r/M_{r-1})$. If $M$ admits a Cohen-Macaulay filtration with respect to $\fa$, then we say that $M$ is \textit{sequentially Cohen-Macaulay with respect to $\fa$}. Note that if $(R, \fm')$ is a local ring, then $M$ is sequentially Cohen-Macaulay if and only if $M$ is sequentially Cohen-Macaulay with respect to $\fm'$.

Let $\Bbbk$ be an algebraically closed field. Theorem \ref{4.6} describes the sequentially Cohen-Macaulayness of $L\otimes_{\Bbbk} N$ with respect to $\mathcal{I}$ in terms of the sequentially Cohen-Macaulayness of $L$ and $N$ with respect to $I$ and $J$, respectively. Moreover, if $L\otimes_{\Bbbk} N$ is sequentially Cohen-Macaulay with respect to $\mathcal{I}$, then we have $\f_{\mathcal{I}}(L\otimes_{\Bbbk} N)=\f_I(L)+\f_J(N)$, provided that $\grade(I, L)>0$ and $\grade(J, N)>0$, as shown in Proposition \ref{4.10}.

In particular, in the local case, $L\otimes_{\Bbbk} N$ is sequentially Cohen-Macaulay as an $A\otimes_\Bbbk B$-module if and only if $L$ and $N$ are sequentially Cohen-Macaulay over $A$ and $B$, respectively. A graded version of this result can be seen in \cite[Theorem 2.11]{STY}. Furthermore, if $L\otimes_{\Bbbk} N$ is sequentially Cohen-Macaulay with $\depth_A L >0$ and $\depth_B N>0$, then we have $\f_{\mathfrak{M}}(L\otimes_{\Bbbk} N)=\f_\fm(L)+\f_\fn(N)$. 

\section{Cohen--Macaulayness}

We recall a consequence of the K\"unneth formula.
\begin{fact}{\em
\label{2.1}
Let $\Bbbk$ be a field, and let $A$ and $B$ be $\Bbbk$-algebras. Let $L$ and $L'$ are $A$-modules such that $L$ is finitely generated, and
 let $N$ and $N'$ are $B$-modules such that $N$ is finitely generated. For each $i\geq 0$, there are $A\otimes_{\Bbbk}B$-module isomorphisms
 \begin{itemize}
\item[{(a)}] $\Ext^i_{A\otimes_{\Bbbk}B}(L\otimes_{\Bbbk}N, L'\otimes_{\Bbbk}N')\cong \underset{s+t=i}\bigoplus \Ext^s_A(L, L')\otimes_{\Bbbk}\Ext^t_B(N, N')$.
\item[{(b)}] $\Tor^{A\otimes_{\Bbbk}B}_i(L\otimes_{\Bbbk}N, L'\otimes_{\Bbbk}N')\cong \underset{s+t=i}\bigoplus \Tor^A_s(L, L')\otimes_{\Bbbk}\Tor^B_t(N, N')$.
\end{itemize}
See for example \cite[Proposition A.1.5]{W}.}
\end{fact}

\begin{lemma}
\label{2.2}
Let $\Bbbk$ be a field, and let $A$ and $B$ be Noetherian $\Bbbk$-algebras such that $A \otimes_\Bbbk B$ is Noetherian. 
Let $I$ be an ideal of $A$ and $J$ an ideal of $B$, and assume that $L$ and $N$ are modules over $A$ and $B$, respectively.
Then,  for each $i\geq 0$, there are $A\otimes_{\Bbbk}B$-module isomorphisms
 \begin{itemize}
\item[{(a)}] $\text{H}_{I\otimes_\Bbbk B}^i(L\otimes_{\Bbbk} N)\cong \text{H}_{I}^i(L)\otimes_\Bbbk N$ for $i\geq 0$.
\item[{(b)}] $\text{H}_{A\otimes_\Bbbk J}^i(L\otimes_{\Bbbk} N)\cong L\otimes_\Bbbk \text{H}_{J}^i(N)$ for $i\geq 0$.
\end{itemize}
\end{lemma}
 \begin{prf}
(a): Consider the canonical flat homomorphism of Noetherian rings $f:A \to A\otimes_{\Bbbk} B$. For every non-negative integer $i$, the independent theorem yields a natural $A\otimes_{\Bbbk} B$-isomorphism
\[
  \text{H}_{I\otimes_\Bbbk B}^i(L\otimes_{\Bbbk} N)  \cong \text{H}_{I}^i(L\otimes_{\Bbbk} N),
\]
via the canonical ring homomorphism $f$.
For each $i\geq 0$, using Fact \ref{2.1}(a), we have the following isomorphisms of $A\otimes_{\Bbbk} B$-modules.
$$\begin{array}{ll}
\text{H}_{I}^i(L\otimes_{\Bbbk} N) & \cong \varinjlim \limits_{n\in \mathbb{N}} \Ext_{A\otimes_{\Bbbk}B}^i\big((A\otimes_{\Bbbk}B)/I^n, L\otimes_{\Bbbk} N\big)\\
                          &\cong  \varinjlim \limits_{n\in \mathbb{N}} \Ext_{A\otimes_{\Bbbk}B}^i(A/I^n\otimes_{\Bbbk}B , L\otimes_{\Bbbk}N)\\
                          &\cong \underset{s+t=i}\bigoplus\varinjlim \limits_{n\in \mathbb{N}}\Ext^s_{A}(A/I^n, L)\otimes_{\Bbbk}  \Ext^t_{B}(B, N)\\
                         & \cong   \text{H}_{I}^i(L)\otimes_\Bbbk N. \\
\end{array}$$
Consequently,
\[
\text{H}_{I\otimes_\Bbbk B}^i(L\otimes_{\Bbbk} N)\cong \text{H}_{I}^i(L)\otimes_\Bbbk N,
\]
for $i\geq 0$, as desired. Part (b) is proved in the same way.
\end{prf}

Let $R$ be a Noetherian ring, $\fa$ an ideal of $R$, and $M$ an $R$-module. The cohomological dimension of $M$ with respect to $\fa$, denoted by $\cd(\fa, M)$, 
is defined as the supremum of
all integers $i$ for which $\text{H}^i_\fa(M)\neq 0$. It is important to note that if $M$ is a finitely generated $R$-module, then $M=\fa M$ if and only if $\cd(\fa, M)=-\infty$.

Now, let $M$ be a finitely generated $R$-module such that $\fa M\neq M$. The grade of $M$ with respect to $\fa$, denoted by $\grade(\fa, M)$,
 is defined as the infimum of
all integers $i$ such that $\text{H}^i_\fa(M)\neq 0$. It is worth mentioning that $M=\fa M$ if and only if $\grade(\fa, M)=+\infty$.
(Recall that, by convension, $\inf \emptyset=+\infty$ and $\sup \emptyset=-\infty$). 
\begin{remark}
\label{2.3}{\em
Let $\Bbbk$ be a field, and let $A$ and $B$ be $\Bbbk$-algebras. Let $I$ be an ideal of $A$ and $J$ an ideal of $B$, and assume that $L$ and $N$ are finitely generated
modules over $A$ and $B$, respectively. We set $\mathcal{I}= I\otimes_\Bbbk B + A \otimes_\Bbbk J$.  
Observe that, $\mathcal{I}(L\otimes_{\Bbbk} N)\neq L\otimes_{\Bbbk} N$ if and only if $IL\neq L$ and $JN\neq N$.
In fact,
$$\begin{array}{ll}
\mathcal{I}(L\otimes_{\Bbbk} N)= L\otimes_{\Bbbk} N & \Longleftrightarrow (L\otimes_{\Bbbk} N)/\mathcal{I}(L\otimes_{\Bbbk} N)=0\\
                                                    &\Longleftrightarrow ((A\otimes_{\Bbbk} B)/\mathcal{I}) \otimes_{A\otimes_{\Bbbk} B} (L\otimes_{\Bbbk} N)= 0 \\
                                                     &\Longleftrightarrow (A/I\otimes_\Bbbk B/J) \otimes_{A\otimes_{\Bbbk} B} (L\otimes_{\Bbbk} N)= 0 \\
                                                      &\Longleftrightarrow  (A/I\otimes_A L) \otimes_\Bbbk (B/J\otimes_B N)=0 \\
                                                    &\Longleftrightarrow  (L/IL) \otimes_\Bbbk (N/JN)=0 \\
                                                    &\Longleftrightarrow L=IL \; or \; N=JN.\\
\end{array}$$
The third step follows from the following fact that
\begin{eqnarray}	
\label{2.3.1}
\frac{A\otimes_{\Bbbk} B}{I \otimes_\Bbbk B+ A\otimes_\Bbbk J}\cong A/I\otimes_\Bbbk B/J.
\end{eqnarray}	
The fourth step follows from Fact \ref{2.1}(b), and the other steps are standard.
}
\end{remark}
From Lemma \ref{2.2}, one obtaines
\begin{corollary}
\label{2.4}
Let $\Bbbk$ be a field, and let $A$ and $B$ be two Noetherian $\Bbbk$-algebras such that $A \otimes_\Bbbk B$ is Noetherian.
 Let $I$ be an ideal of $A$ and $J$ an ideal of $B$, and assume that $L$ and $N$ are nonzero finitely generated
modules over $A$ and $B$, respectively. Then
 \begin{itemize}
\item[{(a)}] $\grade(I \otimes_\Bbbk B, L\otimes_{\Bbbk} N)=\grade(I, L)$  and $\cd(I \otimes_\Bbbk B, L\otimes_{\Bbbk} N)=\cd(I, L)$.
\item[{(b)}] $\grade(A \otimes_\Bbbk J, L\otimes_{\Bbbk} N)=\grade(J, N)$  and  $\cd(A \otimes_\Bbbk J, L\otimes_{\Bbbk} N)=\cd(J, N)$.
\end{itemize}
\end{corollary}
\begin{prf}
(a): We first assume that $( I \otimes_\Bbbk B)(L\otimes_{\Bbbk} N)= L\otimes_{\Bbbk} N$.  This is equivalent to saying that $L=IL$ by Remark \ref{2.3}.
Thus, $\grade(I \otimes_\Bbbk B, L\otimes_{\Bbbk} N)=\grade(I, L)=+\infty$, and so the desired equality holds in this case.  Thus, we may assume that  $( I \otimes_\Bbbk B)(L\otimes_{\Bbbk} N)\neq L\otimes_{\Bbbk} N$. This is equivalent to saying that $L\neq IL$ by Remark \ref{2.3}. Thus, in view of Lemma \ref{2.2}(a), we have 
$$\begin{array}{ll}
\grade(I \otimes_\Bbbk B, L\otimes_{\Bbbk} N)& =\inf\{ i: \text{H}_{I\otimes_\Bbbk B}^i(L\otimes_{\Bbbk} N)\neq 0\}\\
                                                                     &= \inf \{ i: \text{H}_{I}^i(L)\otimes_\Bbbk N\neq 0\}\\
                                                                      &= \inf \{ i: \text{H}_{I}^i(L)\neq 0\}\\
                                                                      &=\grade(I, L).                                                                       
\end{array}$$
Similarly, if $( I \otimes_\Bbbk B)(L\otimes_{\Bbbk} N)= L\otimes_{\Bbbk} N$, then $\cd(I \otimes_\Bbbk B, L\otimes_{\Bbbk} N)=\cd(I, L)=-\infty$, and so the desired equality holds. If 
$( I \otimes_\Bbbk B)(L\otimes_{\Bbbk} N)\neq L\otimes_{\Bbbk} N$, then in view of Lemma \ref{2.2}(b), we have
 $$\begin{array}{ll}
\cd(I \otimes_\Bbbk B, L\otimes_{\Bbbk} N)& =\sup\{ i: \text{H}_{I\otimes_\Bbbk B}^i(L\otimes_{\Bbbk} N)\neq 0\}\\
                                                                     &= \sup \{ i: \text{H}_{I}^i(L)\otimes_\Bbbk N\neq 0\}\\
                                                                      &= \sup \{ i: \text{H}_{I}^i(L)\neq 0\}\\
                                                                      &=\cd(I, L).                                                                       
\end{array}$$
Part (b) is proved similarly. 
\end{prf}
We remark that \cite[Theorem 1.1]{BK}(a) can be deduced from Corollary \ref{2.4}.

Let $M$ be an $R$-module and let $\fa$ and $\fb$ be two finitely generated ideals of $R$. Suppose $\boldsymbol{x}$ and $\boldsymbol{y}$ are two finite sets of generators of $\fa$ and $\fb$, respectively.
Then, the double complex with components $\mathcal{C}^i_{\boldsymbol{x}}(R)\otimes_R\mathcal{C}^j_{\boldsymbol{y}}(M)$ gives rise to a spectral sequence:
\begin{eqnarray}  
\label{2.4.1}
E_2^{i,j}=\text{H}^i_\fa(\text{H}^j_\fb(M))\underset{i}\Longrightarrow \text{H}^{i+j}_{\fa+\fb} (M),
\end{eqnarray}
where the second term $\text{H}^i_\fa(\text{H}^j_\fb(M))$ abuts to $\text{H}^{i+j}_{\fa+\fb} (M)$.

\begin{proposition}
\label{2.5}
Let $\Bbbk$ be a field, and let $A$ and $B$ be two Noetherian $\Bbbk$-algebras such that $A \otimes_\Bbbk B$ is Noetherian.
 Let $I$ be an ideal of $A$ and $J$ an ideal of $B$, and assume that $L$ and $N$ are finitely generated
modules over $A$ and $B$, respectively. Suppose $IL\neq L$ and $JN\neq N$.
Then there are $A\otimes_{\Bbbk}B$-module isomorphisms
\begin{itemize}
\item[{(a)}]$\text{H}_I^{\grade(I, L)}(L)\otimes_\Bbbk \text{H}_{J}^{\grade(J, N)}(N)\cong \text{H}^{\grade(I, L)+\grade(J, N)}_{I \otimes_\Bbbk B + A \otimes_\Bbbk J} (L\otimes_{\Bbbk} N)$.

\item[{(b)}] $\text{H}_I^{\cd(I, L)}(L)\otimes_\Bbbk \text{H}_{J}^{\cd(J, N)}(N)\cong \text{H}^{\cd(I, L)+\cd(J, N)}_{I \otimes_\Bbbk B + A \otimes_\Bbbk J} (L\otimes_{\Bbbk} N)$.
\end{itemize}
\end{proposition}
\begin{prf}
Since $A\otimes_{\Bbbk}B$ is Noetherian, the ideals $I\otimes_\Bbbk B$ and $A\otimes_\Bbbk J$ are finitely generated ideals of $A\otimes_{\Bbbk}B$. In view of (\ref{2.4.1}),
we have the following spectral sequence
\[
E_2^{i,j}= \text{H}^i_{I \otimes_\Bbbk B}\big(\text{H}^j_{A\otimes_\Bbbk J}(L\otimes_{\Bbbk} N)\big)\underset{i}\Longrightarrow \text{H}^{i+j}_{I \otimes_\Bbbk B + A \otimes_\Bbbk J} (L\otimes_{\Bbbk} N).
\]
 Observe that
$$\begin{array}{ll}
\text{H}^i_{I \otimes_\Bbbk B}\big(\text{H}^j_{A\otimes_\Bbbk J}(L\otimes_{\Bbbk} N)\big) & \cong \text{H}^i_{I \otimes_\Bbbk B}\big(L\otimes_\Bbbk \text{H}_{J}^j(N)\big)\\
                                                                                           & \cong \text{H}^i_I(L)\otimes_\Bbbk \text{H}_{J}^j(N).
\end{array}$$
The first step follows from Lemma \ref{2.2}(b) and the second step is by Lemma \ref{2.2}(a).
 We set  $\mathcal{I}= I \otimes_\Bbbk B + A \otimes_\Bbbk J$. Thus we have the following spectral sequence
 \begin{eqnarray}	
\label{2.5.1}
E_2^{i,j}\cong \text{H}^i_I(L)\otimes_\Bbbk \text{H}_{J}^j(N)\underset{i}\Longrightarrow \text{H}^{i+j}_{\mathcal{I}} (L\otimes_{\Bbbk} N).
\end{eqnarray}
We note that for the $E_2$-terms in the spectral sequence we have,
\[
E_2^{i,j}\cong \text{H}^i_I(L)\otimes_\Bbbk \text{H}_{J}^j(N)=0
\]
 if $i<\grade(I, L)$ or  $i>\cd(I, L)$ or $j<\grade(J, N)$ or $j>\cd(J, N)$. 
 Thus the possible nonzero  $E_2$-terms are in the shadowed region of the above picture.
  \begin{figure}
\begin{center}
\psset{unit=1.0cm}
\begin{pspicture}(0,0)(6, 6)
\psline(0.5,1)(0.5, 5)
 \psline(0,1.5)(5.5,1.5)
 \pspolygon[style=fyp, fillcolor=light](1.5,2.)(1.5, 3.9)(4.5, 3.9)(4.5,2.)
 \psline[linestyle=dashed](4.5, 4)(4.5, 2.2 )
 \psline[linestyle=dashed](0.5, 3)(3, 3)
 \psline[linestyle=dashed](3, 3)(3, 1.5)
 \psline[linestyle=dashed](0.0, 4)(0.5, 3.5)
 \psline[linestyle=dashed](0.5,3.5)(1, 3)
 \psline[linestyle=dashed](1, 3)(1.5, 2.5)
 \psline[linestyle=dashed](1.5,2.5)(2, 2)
 \psline[linestyle=dashed](2, 2)(2.5,1.5)
  \psline[linestyle=dashed](2.5, 1.5)(3,1)
 \rput(0.3, 3){\blue{$j$}}
 \rput(3.1,1.3){\blue{$i$}}
 \rput(1.5, 2.){\red{$\bullet$}}
 \rput(4.5, 3.9){\red{$\bullet$}}
 \rput(4.5, 4.2){{$(l, m)$}}
 \rput(3, 3){\red{$\bullet$}}
 \rput(1.4, 1.8){$(s, t)$}
 \rput(3.5, 3){\blue{$E^{i,j}_2$}}
 \end{pspicture}
\end{center}
\end{figure}
Here,  $s=\grade(I, L)$,  $l=\cd(I, L)$, $t=\grade(J, N)$ and $m=\cd(J, N)$.
Considering Figure above,  in order to prove the corner isomorphism (a), let's consider $r\geq 2$, and the following differentials
 \[
 \dots \to E_r^{s-r, t+r-1}\to E_r^{s, t}\to E_r^{s+r, t+1-r}\to\cdots.
 \]
Since $E_2^{p, q}=0$ for $p<s$ or $q<t$, and $E_r^{p, q}$  is a subquotient of $E_2^{p, q}$  for $r\geq 2$, it follows that 
 $E_r^{s-r, t+r-1}=E_r^{s+r, t+1-r}=0$ for $r\geq 2$, and hence $E_2^{s, t}=E_3^{s, t}=\dots=E_\infty^{s, t}.$
 Let's consider the value $n=s+t$.  According to (\ref{2.5.1}), there exists a chain
\[
0=U^{-1}\subseteq U^0\subseteq \dots \subseteq U^n=\text{H}^n_\mathcal{I}(L\otimes_{\Bbbk} N)
\]
of submodules $\text{H}^n_\mathcal{I}(L\otimes_{\Bbbk} N)$ such that $U^p/U^{p-1}\cong E_{\infty}^{p, n-p}$ for $p=0, \dots, n$.
Note that  
\[
E_2^{p, n-p}\cong \text{H}^p_I(L)\otimes_\Bbbk \text{H}_{J}^{n-p}(N)=0 \quad\text{for}\quad p\neq s.
\]
Since $E_{\infty}^{p, n-p}$ is a subquotient of $E_{2}^{p, n-p}$, it follows that $U^0=U^1=\dots=U^{s-1}=0$ and $U^s=U^{s+1}=\dots=U^n$. Consequently,
\[
E_2^{s, t}=E_{\infty}^{s, t}\cong U^s/U^{s-1}\cong U^s=U^n=\text{H}^n_\mathcal{I}(L\otimes_{\Bbbk} N).
\]
This yields the first isomorphism in (a).   The second corner isomorphism (b) can be proven using the same method. 
\end{prf}
As the main result of this section, we extend the result stated in \cite[Theorem 1.1]{BK}(b) and Corollary \ref{2.4} as follows:
\begin{theorem}
\label{2.6}
Let $\Bbbk$ be a field, and let $A$ and $B$ be two Noetherian $\Bbbk$-algebras such that $A \otimes_\Bbbk B$ is Noetherian.
 Let $I$ be an ideal of $A$ and $J$ an ideal of $B$, and assume that $L$ and $N$ are nonzero finitely generated
modules over $A$ and $B$, respectively. Then
 \begin{itemize}
\item[{(a)}] $\grade(I \otimes_\Bbbk B+ A\otimes_\Bbbk J, L\otimes_{\Bbbk} N)=\grade(I, L)+\grade(J, N)$.
\item[{(b)}]  $\cd(I \otimes_\Bbbk B+ A\otimes_\Bbbk J, L\otimes_{\Bbbk} N)=\cd(I, L)+\cd(J, N)$.
\end{itemize}
\end{theorem}
\begin{prf}
 (a):  We define $\mathcal{I}= I \otimes_\Bbbk B + A \otimes_\Bbbk J$. If $\mathcal{I}(L\otimes_{\Bbbk} N)= L\otimes_{\Bbbk} N$, then
 $\grade(\mathcal{I}, L\otimes_{\Bbbk} N)=+\infty$. According to Remark \ref{2.3}, we have $L=IL$ or $N=JN$.  Thus,  $\grade(I, L)=+\infty$ or
 $\grade(J, N)=+\infty$, and so the desired equality holds. 
 Therefore, we can assume that $\mathcal{I}(L\otimes_{\Bbbk} N)\neq L\otimes_{\Bbbk} N$, which implies that $IL\neq L$ and $JN\neq N$.
 We consider the spectral sequence of (\ref{2.5.1})
\[
E_2^{i,j}\cong \text{H}^i_I(L)\otimes_\Bbbk \text{H}_{J}^j(N)\underset{i}\Longrightarrow \text{H}^{i+j}_{\mathcal{I}} (L\otimes_{\Bbbk} N),
\]
and set  $t=\min \{i+j: E_2^{i,j}\neq 0\}$. We claim that $\grade(\mathcal{I}, L\otimes_{\Bbbk} N) = t$. First, we show that $\grade(\mathcal{I}, L\otimes_{\Bbbk} N) \geq  t$.
 In order to show this, we claim that  $\text{H}^n_\mathcal{I}(L\otimes_{\Bbbk} N)=0$ for $n< t$.
 For the non-negative integer $n$, there exists a chain
\[
0=U^{-1}\subseteq U^0\subseteq \dots \subseteq U^n=\text{H}^n_\mathcal{I}(L\otimes_{\Bbbk} N)
\]
of submodules $\text{H}^n_\mathcal{I}(L\otimes_{\Bbbk} N)$ such that $U^p/U^{p-1}\cong E_{\infty}^{p, n-p}$ for $p=0, \dots, n$.
Since $n< t$, we have $E_2^{p, n-p}=0$ for $p=0, \dots, n$. Moreover, $E_{\infty}^{p, n-p}$ is a subquotient of $E_{2}^{p, n-p}$, so it follows that $E_{\infty}^{p, n-p}=0$ for $p=0, \dots, n$.
Thus, from the chain above, we can deduce that $\text{H}^n_\mathcal{I}(L\otimes_{\Bbbk} N)=0$ for $n< t$.  Therefore, $\grade(\mathcal{I}, L\otimes_{\Bbbk} N) \geq  t$.
We observe that 
$$
\begin{array}{ll}
\grade(\mathcal{I}, L\otimes_{\Bbbk} N)  & \geq  t=\min \{i+j: E_2^{i,j}\neq 0\}\\
                                         & =  \min \{i+j: \text{H}^i_I(L)\otimes_\Bbbk \text{H}_{J}^j(N)\neq 0\}\\
                                         & =  \min \{i: \text{H}^i_I(L)\neq 0\}+ \min \{j: \text{H}_{J}^j(N)\neq 0\}\\
                                         & =\grade(I, L)+\grade(J, N).\\
\end{array}
$$
 On the other hand, as $\text{H}_I^{\grade(I, L)}(L)\neq 0$ and $\text{H}_{J}^{\grade(J, N)}(N)\neq 0$, it follows from Proposition \ref{2.5}(a) that
  \[
  \text{H}^{\grade(I, L)+\grade(J, N)}_{\mathcal{I}} (L\otimes_{\Bbbk} N)\neq 0.
  \]
 Consequently,
 \[
 \grade(\mathcal{I}, L\otimes_{\Bbbk} N)\leq \grade(I, L)+\grade(J, N).
 \]
  Therefore, the equality (a) holds.

(b): Similar to the first part, if $\mathcal{I}(L\otimes_{\Bbbk} N)= L\otimes_{\Bbbk} N$, then
 $\cd(\mathcal{I}, L\otimes_{\Bbbk} N)=-\infty$. By Remark \ref{2.3}, we have $L=IL$ or $N=JN$.  Hence,  $\cd(I, L)=-\infty$ or $\cd(J, N)=-\infty$
 and so the desired equality holds. Thus, we may assume that $\mathcal{I}(L\otimes_{\Bbbk} N)\neq L\otimes_{\Bbbk} N$,  and so $IL\neq L$ and $JN\neq N$.
 Let $t=\max \{i+j: E_2^{i,j}\neq 0\}$. We claim that $\cd(\mathcal{I}, L\otimes_{\Bbbk} N)  =t$. First, we show that $\cd(\mathcal{I}, L\otimes_{\Bbbk} N)  \leq t$. In order to show this,
we claim that $\text{H}^n_\mathcal{I}(L\otimes_{\Bbbk} N)=0$ for $n > t$. According to (\ref{2.5.1}), there exists a chain
\[
0=U^{-1}\subseteq U^0\subseteq \dots \subseteq U^n=\text{H}^n_\mathcal{I}(L\otimes_{\Bbbk} N)
\]
of submodules $\text{H}^n_\mathcal{I}(L\otimes_{\Bbbk} N)$ such that $U^p/U^{p-1}\cong E_{\infty}^{p, n-p}$ for $p=0, \dots, n$.
Since $n> t$, we have $E_2^{p, n-p}=0$ for $p=0, \dots, n$. Moreover, $E_{\infty}^{p, n-p}$ is a subquotient of $E_{2}^{p, n-p}$, so it follows that $E_{\infty}^{p, n-p}=0$ for $p=0, \dots, n$.
Thus, from the chain above, we can deduce that $\text{H}^n_\mathcal{I}(L\otimes_{\Bbbk} N)=0$ for $n> t$.
Therefore, $\cd(\mathcal{I}, L\otimes_{\Bbbk} N)  \leq t$. Observe that
$$\begin{array}{ll}
\cd(\mathcal{I}, L\otimes_{\Bbbk} N)  & \leq t=\max \{i+j: E_2^{i,j}\neq 0\}\\
                                      & =  \max \{i+j: \text{H}^i_I(L)\otimes_\Bbbk \text{H}_{J}^j(N)\neq 0\}\\
                                      & =  \max \{i: \text{H}^i_I(L)\neq 0\}+ \max \{j: \text{H}_{J}^j(N)\neq 0\}\\
                                      & =\cd(I, L)+\cd(J, N).\\
\end{array}$$
On the other hand, since $\text{H}_I^{\cd(I, L)}(L)\neq 0$ and $\text{H}_{J}^{\cd(J, N)}(N)\neq 0$, it follows from Proposition \ref{2.5}(b) that
\[
\text{H}^{\cd(I, L)+\cd(J, N)}_{\mathcal{I}} (L\otimes_{\Bbbk} N)\neq 0.
\]
Consequently, $\cd(\mathcal{I}, L\otimes_{\Bbbk} N)\geq \cd(I, L)+\cd(J, N)$. Therefore, the equality (b) holds.
\end{prf}
 Let $R$ be a Noetherian ring, $\fa$ a proper ideal of $R$, and $M$ a finitely generated $R$-module. We say that $M$ is \textit{Cohen--Macaulay with respect to $\fa$} if either $M=\fa M$ or $M\neq \fa M$ and $\grade(\fa, M)=\cd(\fa, M)$, as defined in \cite{R2}. It is important to note that when $(R, \fm')$ is a local ring, $M$ is Cohen-Macaulay if and only if $M$ is Cohen-Macaulay with respect to $\fm'$.

As a consequence of Theorem \ref{2.6}, we have
\begin{corollary}
\label{2.7}
Suppose $IL\neq L$ and $JN\neq N$. Then  $L\otimes_{\Bbbk} N$ is Cohen--Macaulay with respect to $I \otimes_\Bbbk B+ A\otimes_\Bbbk J$ if and only if
 $L$  and $N$ are Cohen--Macaulay with respect to $I$ and $J$, respectively.
\end{corollary}

\begin{remark}
\label{2.8}{\em
Let $(A, \fm)$ and $(B, \fn)$ be two local algebras over a field $\Bbbk$, and sharing $\Bbbk$ as the common residue field. We also assume that either $A$ or $B$ is algebraic over $\Bbbk$. Then $A\otimes_\Bbbk B$ is a local ring with the unique maximal ideal $\mathfrak{M}= \fm \otimes_\Bbbk B + A \otimes_\Bbbk \fn$. In fact, by (\ref{2.3.1}), we have  
$$\begin{array}{ll}
\frac{A\otimes_{\Bbbk} B}{\fm \otimes_\Bbbk B+ A\otimes_\Bbbk \fn} & \cong A/\fm\otimes_\Bbbk B/\fn\\
                    & \cong  \Bbbk \otimes_{\Bbbk} \Bbbk \\
                    & \cong \Bbbk.\\
\end{array}$$
We set $\mathfrak{M}= \fm \otimes_\Bbbk B + A \otimes_\Bbbk \fn$. Thus,  $\mathfrak{M}$  is a maximal ideal of $A\otimes_\Bbbk B$.  Moreover, $A/\fm\otimes_\Bbbk B/\fn$ is a local ring. Hence, \cite[Theorem 4]{L} yields that $A\otimes_\Bbbk B$ is a local ring; thus $\mathfrak{M}$ is the unique maximal ideal of $A\otimes_\Bbbk B$.
}
\end{remark}
As a result of Theorem \ref{2.6}, we also obtain the following consequence. A graded version of this fact is presented in \cite[Corollary 2.3]{STY}.
\begin{corollary}
\label{2.9}
Let $(A, \fm)$ and $(B, \fn)$ be two Noetherian local algebras over a field $\Bbbk$ such that $A \otimes_\Bbbk B$ is Noetherian.  We also suppose 
that $\Bbbk$ is the common residue field of both $A$ and $B$, and either $A$ or $B$ is algebraic over $\Bbbk$.   Let $L$ and $N$ be 
two nonzero finitely generated modules over  $A$ and $B$, respectively. Then
 \begin{itemize}
\item[{(a)}] $\depth_{A \otimes_\Bbbk B} L\otimes_{\Bbbk} N=\depth_A L+\depth_B N$.
\item[{(b)}] $\dim_{A \otimes_\Bbbk B} L\otimes_{\Bbbk} N=\dim_A L+\dim_B N$.
\end{itemize}
\end{corollary}
\begin{prf}
The assertions follow from Remark \ref{2.8} and Theorem \ref{2.6}.
\end{prf}

\section{The finiteness dimension}
Let $R$ be a Noetherian ring, $\fa$ an ideal of $R$, and $M$ a finitely generated $R$-module. The \emph{finiteness dimension of $M$ with respect to $\fa$}, denoted by $\f_\fa(M)$, is defined as the infimum of non-negative integers $i$ for which $\text{H}_\fa^i(M)$ is not finitely generated. Note that $\f_\fa(M)$ can be either a positive integer or $+\infty$. 
If $\cd(\fa, M)\geq 1$, then $\text{H}^{\cd(\fa, M)}_\fa(M)$ is not finitely generated, see for instance \cite[Corollary 3.3(i)]{DV}.
Thus, in this case, one has $\f_\fa(M)\leq \cd(\fa, M)$.
\begin{definition}
\label{3.1}
We say $M$ is {\it generalized Cohen--Macaulay with respect to $\fa$} if $\cd(\fa, M)\leq 0$ or $\f_\fa(M)=\cd(\fa, M)$, as defined in \cite{DGTZ}.
\end{definition}
We observe that if $(R, \fm')$ is a local ring, then $M$ is generalized Cohen--Macaulay if and only if $M$ is generalized Cohen--Macaulay with respect to $\fm'$.

Let $\Bbbk$ be a filed, and let $A$ and $B$ be two Noetherian local $\Bbbk$-algebras such that $A \otimes_\Bbbk B$ is Noetherian.
  Let $I$ be an ideal of $A$ and $J$ an ideal of $B$, and assume that $L$ and $N$ are nonzero finitely generated
modules over $A$ and $B$, respectively.  From Lemma \ref{2.2} one obtaines 
\[
\f_{I \otimes_\Bbbk B}(L\otimes_{\Bbbk} N)=\f_I(L) \;\; \text {and} \;\; \f_{A \otimes_\Bbbk J}(L\otimes_{\Bbbk} N)=\f_J(N).
\]
In fact, 
$$\begin{array}{ll}
\f_{I \otimes_\Bbbk B}(L\otimes_{\Bbbk} N)& = \inf \{ i\in \mathbb{N}_0: \text{H}^i_{I \otimes_\Bbbk B}(L\otimes_{\Bbbk} N) \; \text {is not finitely generated} \; \}\\
&=\inf \{ i\in \mathbb{N}_0: \text{H}_{I}^i(L)\otimes_\Bbbk N \; \text {is not finitely generated} \; \}\\
&=\inf \{ i\in \mathbb{N}_0: \text{H}_{I}^i(L) \; \text {is not finitely generated} \; \}\\
&= \f_{I} (L).\\
\end{array}$$
The second part can be seen similarly. 
As the main result of this section, we extend this fact as follows:
\begin{theorem}
\label{3.2}
Let $\Bbbk$ be a filed, and let $A$ and $B$ be two Noetherian local $\Bbbk$-algebras such that $A \otimes_\Bbbk B$ is Noetherian.
  Let $I$ be an ideal of $A$ and $J$ an ideal of $B$, and assume that $L$ and $N$ are nonzero finitely generated
modules over $A$ and $B$, respectively.  Then
\[
\f_{I \otimes_\Bbbk B+A \otimes_\Bbbk J}(L\otimes_{\Bbbk} N)= \min\{\grade(I, L)+\f_J(N), \f_I(L)+\grade(J, N)\}.
\]
\end{theorem}
\begin{prf}
We define $\mathcal{I}=I \otimes_\Bbbk B+A \otimes_\Bbbk J$ and initially assume that $\mathcal{I}(L\otimes_{\Bbbk} N) = L\otimes_{\Bbbk} N$. Thus, $\text{H}^i_\mathcal{I}(L\otimes_{\Bbbk} N)=0$ for all $i$, implying that  $\f_{\mathcal{I}}(L\otimes_{\Bbbk} N)=+\infty$. As per Remark \ref{2.3}, we have $IL= L$ or $JN=N$. If $IL= L$ and $JN=N$, then $\grade(I, L)=\grade(J, N)=+\infty$ and $\f_I(L)=\f_J(N)=+\infty$. Thus the desired equality holds in this case. If $IL= L$ and $JN\neq N$, then $\grade(I, L)=\f_I(L)=+\infty$,  $\grade(J, N)<\infty$, and $\f_J(N)$ is a positive integer or $+\infty$. Thus, the required equality holds in this case as well. 
Therefore, we can assume that $\mathcal{I}(L\otimes_{\Bbbk} N)\neq L\otimes_{\Bbbk} N$, and so $\cd(\mathcal{I}, L\otimes_{\Bbbk} N)\geq 0$. If $\cd(\mathcal{I}, L\otimes_{\Bbbk} N)=0$, then by Theorem \ref{2.6}(b), we have $\cd(I, L)=\cd(J, N)=0$. Consequently, $\f_{\mathcal{I}}(L\otimes_{\Bbbk} N)=\f_I(L)=\f_J(N)=+\infty$, and therefore the claimed equality holds in this case as well. Thus, we may assume that $\cd(\mathcal{I}, L\otimes_{\Bbbk} N)\geq 1$.
Consider the spectral sequence of (\ref{2.5.1})
\[
E_2^{i,j}\cong \text{H}^i_I(L)\otimes_\Bbbk \text{H}_{J}^j(N)\underset{i}\Longrightarrow \text{H}^{i+j}_{\mathcal{I}} (L\otimes_{\Bbbk} N),
\]
and set $t=\min \{i+j: E_2^{i,j} \; \text {is not finitely generated} \}$. We claim that $f_\mathcal{I}(L\otimes_{\Bbbk} N)= t$.
First, we show that $f_\mathcal{I}(L\otimes_{\Bbbk} N)\geq t$.
 In order to show this, we claim that $\text{H}^n_\mathcal{I}(L\otimes_{\Bbbk} N)$ is finitely generated for $n < t$.  In light of (\ref{2.5.1}), a chain
\[
0=U^{-1}\subseteq U^0\subseteq \dots \subseteq U^n=\text{H}^n_\mathcal{I}(L\otimes_{\Bbbk} N)
\]
of submodules $\text{H}^n_\mathcal{I}(L\otimes_{\Bbbk} N)$ exists such that $U^p/U^{p-1}\cong E_{\infty}^{p, n-p}$ for $p=0, \dots, n$.
Since $n< t$, it follows that $E_2^{p, n-p}$ is finitely generated for $p=0, \dots, n$.
 Moreover, $E_{\infty}^{p, n-p}$ is a subquotient of $E_{2}^{p, n-p}$, so it follows that $E_{\infty}^{p, n-p}$ is also finitely generated for $p=0, \dots, n$.
 Let $p=0$. Thus, $U^0\cong U^0/U^{-1}=E_{\infty}^{0, n}$ is finitely generated. For $p=1$, we have $U^1/U^{0}\cong E_{\infty}^{1, n-1}$  is finitely generated. So,
 the exact sequence $0\to U^0\to U^1 \to U^1/U^0 \to 0$ yields $U^1$ is finitely generated. Continuing in this way, we obtain that $U^n=\text{H}^n_\mathcal{I}(L\otimes_{\Bbbk} N)$ is
 finitely generated for $n<t$. Therefore, $\f_{\mathcal{I}}(L\otimes_{\Bbbk} N)\geq t$.
 
 In order to show the reverse inequality,  $f_\mathcal{I}(L\otimes_{\Bbbk} N)\leq t$, it suffices to show that
 $\text{H}^t_\mathcal{I}(L\otimes_{\Bbbk} N)$ is not finitely generated.  According to (\ref{2.5.1}), there exists a chain
\[
0=U^{-1}\subseteq U^0\subseteq \dots \subseteq U^t=\text{H}^t_\mathcal{I}(L\otimes_{\Bbbk} N)
\]
of submodules $\text{H}^t_\mathcal{I}(L\otimes_{\Bbbk} N)$ such that $U^p/U^{p-1}\cong E_{\infty}^{p, t-p}$ for $p=0, \dots, t$.  For $p=t$, we have  the following exact sequence, 
 \[
 0\to U^{t-1}\to \text{H}^t_\mathcal{I}(L\otimes_{\Bbbk} N) \to E_{\infty}^{t, 0}\to 0. 
 \]
To complete our proof, we only need to show that $E_{\infty}^{t, 0}$ is not finitely generated.
Let's consider $r\geq 2$ and the following differentials
 \[
 \dots \to E_r^{t-2r,2r-2}\xrightarrow{d_r^{t-2r,2r-2}} E_r^{t-r, r-1}\xrightarrow{d_r^{t-r, r-1}} E_r^{t, 0}\xrightarrow{d_r^{t, 0}} E_r^{t+r, 1-r}=0.
 \]
Here,  $E_r^{t+r, 1-r}=0$  for $r\geq 2$ because $E_r^{t+r, 1-r}$ is a subquotient of $E_2^{t+r, 1-r}$, and $E_2^{t+r, 1-r}=0$  for $r\geq 2$.
As 
\[
E_{r+1}^{t, 0}=\ker d_r^{t, 0}/\im d_r^{t-r, r-1}= E_r^{t, 0}/\im d_r^{t-r, r-1},
\]
we have the following exact sequence
 \begin{eqnarray}	
\label{3.2.1}
 0\to \im d_r^{t-r, r-1}\to E_r^{t, 0} \to E_{r+1}^{t, 0}\to 0.
 \end{eqnarray}
Note that if $i+j<t$, then $E_2^{i, j}$ is finitely generated. Hence $E_r^{i, j}$ is finitely generated for $i+j<t$ and $r\geq 2$. 
In particular, $E_r^{t-r, r-1}$ is finitely generated, and so  $\im d_r^{t-r, r-1}$ is finitely generated. 
  Now, on the contrary suppose $E_{\infty}^{t, 0}$ is finitely generated. Since   $E_{r}^{t, 0}=E_{\infty}^{t, 0}$ for $r$  sufficiently large, 
 we have that  $E_{r}^{t, 0}$ is finitely generated. Fix $r$ and suppose  $E_{r+1}^{t, 0}$ is finitely generated. From the exact sequence (\ref{3.2.1}), we obtain that $E_{r}^{t, 0}$ is finitely generated. 
 Continuing in this fashion, we see that $E_{r}^{t, 0}$ is finitely generated for all $r\geq 2$.  In particular, $E_{2}^{t, 0}$ is finitely generated, a contradiction.  Consequently,
 $\f_{\mathcal{I}}(L\otimes_{\Bbbk} N)\leq t$. Notice that 
 $$\begin{array}{ll}
   t & =\min \{i+j: E_2^{i,j} \; \text {is not finitely generated} \}\\
                                        & =  \min \{i+j: \text{H}^i_I(L)\otimes_\Bbbk \text{H}_{J}^j(N)\; \text {is not finitely generated} \}\\
                                         & =  \min \big \{ \min \{i: \text{H}^i_I(L)\neq 0 \}+ \min \{j: \text{H}_{J}^j(N) \; \text {is not finitely generated}\}, \\
                                        &  \min \{i: \text{H}^i_I(L)\; \text {is not finitely generated} \}+ \min \{j: \text{H}_{J}^j(N)\neq 0 \} \big\}\\
                                         & =\min\{\grade(I, L)+\f_J(N), \f_I(L)+\grade(J, N)\}.\\
\end{array}$$
Therefore, the proof is complete. 
\end{prf}

\begin{corollary}
\label{3.3}
Let $(A, \fm)$ and $(B, \fn)$ be two Noetherian local algebras over a field $\Bbbk$ such that $A \otimes_\Bbbk B$ is Noetherian.  We also suppose 
that $\Bbbk$ is the common residue field of both $A$ and $B$,  and either $A$ or $B$ is algebraic over $\Bbbk$.  Let $L$ and $N$ are nonzero finitely generated
modules over $A$ and $B$, respectively. Then 
\[
\f_{\mathfrak{M}}(L\otimes_{\Bbbk} N)= \min\{\depth_A L+\f_\fn(N), \f_\fm (L)+\depth_B N \}.
\]
\end{corollary}
\begin{prf}
The assertion follows from Remark \ref{2.8} and Theorem \ref{3.2}.
\end{prf}
\begin{proposition}
\label{3.4}
Let $\Bbbk$ be a field , and let $A$ and $B$ be $\Bbbk$-algebras. Let $I$ be an ideal of $A$ and $J$ an ideal of $B$, and assume that $L$ and $N$ are nonzero finitely generated
modules over $A$ and $B$, respectively. Suppose  $\cd(I, L)\geq 1$ and $\cd(J, N)\geq 1$. Then the following conditions are equivalent:
 \begin{itemize}
\item[{(a)}] $L\otimes_{\Bbbk} N$ is generalized Cohen--Macaulay with respect to $I \otimes_\Bbbk B+ A\otimes_\Bbbk J$;
\item[{(b)}] $L\otimes_{\Bbbk} N$ is Cohen--Macaulay with respect to $I \otimes_\Bbbk B+ A\otimes_\Bbbk J$;
\item[{(c)}] $L$ and $N$ are Cohen--Macaulay with respect to $I$ and $J$, respectively.
\end{itemize}
Moreover, if one of these equivalent conditions holds, then we have
\[
\f_{\mathcal{I}}(L\otimes_{\Bbbk} N)=\f_I(L)+\f_J(N),
\]
where  $\mathcal{I}=I \otimes_\Bbbk B+A \otimes_\Bbbk J$.
\end{proposition}
\begin{prf}
$(a)\Rightarrow (c)$: Suppose $L\otimes_{\Bbbk} N$ is generalized
Cohen--Macaulay with respect to $\mathcal{I}$.  As  $\cd(I, L)\geq 1$ and  $\cd(J, N)\geq 1$, we have
$\cd(\mathcal{I}, L\otimes_{\Bbbk} N)\geq 2$. Observe that
$$\begin{array}{ll}
\cd(I, L)+\cd(J, N) & = \cd(\mathcal{I}, L\otimes_{\Bbbk} N)\\
                    & =  \f_{\mathcal{I}}(L\otimes_{\Bbbk} N)\\
                    & = \min\{\grade(I, L)+\f_J(N), \f_I(L)+\grade(J, N)\}\\
                    & \leq \grade(I, L)+\f_J(N)\\
                    & \leq \cd(I, L)+\cd(J, N).\\
\end{array}$$
The first step is by Theorem \ref{2.6}(b). Our assumption yields the second step. Since,  $\cd(I, L)\geq 1$ and  $\cd(J, N)\geq 1$, it follows that 
$\grade(I, L),  \grade(J, N),  \f_I(L) \;\;  \text{and} \;\;  \f_J(N) $ are finite numbers. Hence the third step is by Theorem \ref{3.2}, and the remaining steps are standard.
Consequently,
\[
\cd(I, L)+\cd(J, N)=\grade(I, L)+\f_J(N).
\]
Similarly, one obtains
\[
\cd(I, L)+\cd(J, N)=\f_I(L)+\grade(J, N).
\]
Therefore, $\grade(I, L)=\cd(I, L)$ and $\grade(J, N)=\cd(J, N)$;  so that $L$ and $N$ are Cohen--Macaulay with respect to $I$ and $J$, respectively,  as desired.

$(c)\Longleftrightarrow (b)$: follows from Corollary \ref{2.7}.

$(b)\Rightarrow (a)$:  is clear.
\end{prf}
A graded version of the following result can be seen in \cite[Theorem 2.6]{STY}. 
\begin{corollary}
\label{3.5}
Let $(A, \fm)$ and $(B, \fn)$ be two Noetherian local algebras over a field $\Bbbk$ such that $A \otimes_\Bbbk B$ is Noetherian.  We also suppose 
that $\Bbbk$ is the common residue field of both $A$ and $B$, and either $A$ or $B$ is algebraic over $\Bbbk$.   Let $L$ and $N$ be nonzero finitely generated
modules over $A$ and $B$, respectively. Assume that both $\dim_A L$ and $\dim_B N$ are positive. Then the following conditions are equivalent:
 \begin{itemize}
\item[{(a)}] $L\otimes_{\Bbbk} N$ is generalized Cohen--Macaulay;
\item[{(b)}] $L\otimes_{\Bbbk} N$ is Cohen--Macaulay;
\item[{(c)}] $L$ and $N$ are Cohen--Macaulay.
\end{itemize}
Moreover, if one of these equivalent conditions holds, then we have
\[
\f_{\mathfrak{M}}(L\otimes_{\Bbbk} N)=\f_\fm(L)+\f_\fn(N),
\]
where $\mathfrak{M} = \fm \otimes_\Bbbk B + A \otimes_\Bbbk \fn$.
\end{corollary}
\begin{prf}
The assertions follow from Remark \ref{2.8} and Proposition \ref{3.4}.
\end{prf}
\begin{remark}
\label{3.6}{\em
We observed that, if one of the equivalent conditions of Proposition \ref{3.4} holds, then we have
\[
\f_{\mathcal{I}}(L\otimes_{\Bbbk} N)=\f_I(L)+\f_J(N).
\]
However, the converse is not true in general. In fact, this equality holds for sequentially Cohen--Macaulay modules which is a generalization of Cohen--Macaulay modules, as shown in Proposition \ref{4.10}.
}
\end{remark}
The following example shows that if any of the equivalent conditions stated in Corollary \ref{3.5} does not hold, then 
the equality $\f_{\mathfrak{M}}(L\otimes_{\Bbbk} N)=\f_\fm(L)+\f_\fn(N)$ is no longer valid.
\begin{example}
\label{3.7}{\em
Let $\Bbbk$ be a field. We set $L=\Bbbk[x_1, x_2]/\langle x_1^2, x_1x_2 \rangle$ and consider it as an $\Bbbk[x_1, x_2]$-module. We also set $N=\Bbbk[y_1, y_2]/\langle y_1 \rangle $ and consider it as an $\Bbbk[y_1, y_2]$-module. It can be observed that $L$ has depth $0$ and dimension $1$. Thus $L$ is not Cohen--Macaulay, and  $f_\fm(L)=1$, where $\fm=\langle x_1, x_2 \rangle$ is the unique graded maximal ideal of $\Bbbk[x_1, x_2]$. Moreover, $N$ is Cohen--Macaulay of dimension $1$, so $f_\fn(N)=1$, where $\fn=\langle y_1, y_2 \rangle$ is the unique graded maximal ideal of $\Bbbk[y_1, y_2]$. Now, we consider $L\otimes_\Bbbk N= \Bbbk[x_1, x_2, y_1, y_2]/\langle x_1^2, x_1x_2, y_1\rangle $ as an $\Bbbk[x_1, x_2, y_1, y_2]$-module.  According to Corollary \ref{3.3}, we have
$$\begin{array}{ll}
f_\mathfrak{M}(L\otimes_\Bbbk N) &=\min\{\depth_A L+\f_\fn(N), \f_\fm (L)+\depth_B N \}\\
                                 &=\min\{0+1, 1+1\} \\
                                 &=1,\\
\end{array}$$
where $\mathfrak{M}=\langle x_1, x_2, y_1, y_2\rangle=\fm \otimes_\Bbbk B + A \otimes_\Bbbk \fn$ is the unique graded maximal ideal of $\Bbbk[x_1, x_2, y_1, y_2]$. Here, $A=\Bbbk[x_1, x_2]$ and $B=\Bbbk[y_1, y_2]$. Consequently,
\[
1=f_\mathfrak{M}(L\otimes_\Bbbk N)\neq f_\fm(L)+f_\fn(L)=2.
\]
Moreover,  $L\otimes_\Bbbk N$ is not generalized Cohen--Macaulay, because $ 1=f_\mathfrak{M}(L\otimes_\Bbbk N)\neq \dim_{A\otimes_\Bbbk B} L\otimes_\Bbbk N=2$. }
\end{example}
\section{Sequentially Cohen--Macaulayness}
Let $R$ be a Noetherian ring, $\fa$ an ideal of $R$ and $M$ a finitely generated $R$-module.
\begin{definition}
\label{4.1} {\em A filtration $\mathscr{D}$: $0=D_0\varsubsetneq D_1\varsubsetneq \dots  \varsubsetneq  D_r=M$ of $M$ by submodules $M_i$ is called the
 {\em dimension filtration of $M$ with respect to $\fa$} if $D_{i-1}$ is the largest submodule of $D_i$ for which $\cd(\fa, D_{i-1})<\cd(\fa, D_i)$
for all $i=1, \dots, r$.}
\end{definition}

A finite filtration $\mathcal{F}: 0=M_0\varsubsetneq M_1 \varsubsetneq \dots  \varsubsetneq M_r=M$ of $M$ by submodules $M_i$ is called a
{\em Cohen-Macaulay filtration with respect to $\fa$} if each quotient $M_i/M_{i-1}$ is Cohen-Macaulay with respect to $\fa$ and
$0 \leq \cd(\fa, M_1/M_0)<\cd(\fa, M_2/M_1)< \dots< \cd(\fa, M_r/M_{r-1})$. If $M$ admits a Cohen-Macaulay filtration with respect to $\fa$,
then we say $M$ is \emph {sequentially Cohen-Macaulay with respect to $\fa$}. Note that if $M$ is sequentially Cohen-Macaulay with respect to $\fa$,
then the filtration $\mathscr{F}$ is uniquely determined and it is just the dimension filtration of $M$ with respect to $\fa$,
 i.e., $\mathscr{F}=\mathscr{D}$. The proof of this fact is based on the same method as in the proof of \cite[Proposition 2.9]{R3}
 by replacing the ring $R$ and the general ideal $\fa$ with the polynomial ring $S$ and the ideal $Q$, respectively. 
 Notice that every Cohen--Macaulay module with respect to $\fa$ is sequentially Cohen--Macaulay with respect to $\fa$. Therefore, we may consider
  the zero module to be sequentially Cohen-Macaulay with respect to $\fa$. 
 We observe that if $(R, \fm')$ is a local ring, then $M$ is sequentially Cohen-Macaulay if and only if $M$ is sequentially Cohen-Macaulay with respect to $\fm'$. We recall the following fact which will be used in the sequel.
\begin{fact}{\em
\label{4.2}
Let $\mathscr{D}$ be the dimension filtration of $M$ with respect to $\fa$. Then
\[
\Ass_R D_i/D_{i-1} = \{\fp \in \Ass_R M: \cd(\fa, R/\fp)=\cd(\fa, D_i)\}.
\]
This fact is also proved in the same way as the proof of \cite[Lemma 1.5]{R1} by replacing the ring $R$ and the general ideal $\fa$ with the polynomial ring $S$ and the ideal $Q$, respectively.}
\end{fact}

\begin{fact}{\em
\label{4.3}
The exact sequence  $ 0\rightarrow M' \rightarrow M \rightarrow M'' \rightarrow 0$ of finitely generated $R$-modules yields $ \cd(\fa, M)=\max\{\cd(\fa, M'), \cd(\fa, M'')\}$,
see \cite[Proposition 4.4]{CJR}.}
\end{fact}

Going forward, we will assume that $\Bbbk$ is a field, and that $A$ and $B$ are two Noetherian $\Bbbk$-algebras such that $A \otimes_\Bbbk B$ is Noetherian. Let $I$ and $J$ be two ideals of $A$ and $B$, respectively. As before, we define $\mathcal{I}$ to be $I \otimes_\Bbbk B+ A\otimes_\Bbbk J$. We also assume that  $(A, \fm)$,  $(B, \fn)$  and $\mathfrak{M}= \fm \otimes_\Bbbk B + A \otimes_\Bbbk \fn$ are as described in Remark \ref{2.8}. Additionally, let $L$ and $N$ be nonzero finitely generated modules over $A$ and $B$, respectively.   We shall need the following fact.
\begin{fact}	
\label{4.4}{\em  Continue with the notation and assumptions as above. Let $\Bbbk$ be an algebraically closed field. Then
\[
\Ass_{A \otimes_\Bbbk B} L\otimes_\Bbbk N=\{\fp_1 \otimes_\Bbbk B+ A\otimes_\Bbbk \fp_2: \fp_1 \in \Ass_A  L \quad\text{and}\quad\fp_2\in \Ass_B N\},
\]
see \cite[Corollary 2.8]{HNTT}.}
\end{fact}

 \begin{lemma}
\label{4.5}
Continue with the notation and assumptions as above. Let $\Bbbk$ be an algebraically closed field. Let  $\mathscr{F}$: $0=L_0\varsubsetneq L_1 \varsubsetneq \dots  \varsubsetneq L_r=L$
 be the dimension filtration of $L$ with respect to $I$, and let
$\mathscr{G}$: $0=N_0\varsubsetneq N_1 \varsubsetneq \dots  \varsubsetneq N_s=N$ be the dimension filtration of $N$ with respect to $J$.
Then $\mathscr{G}\otimes_\Bbbk \mathscr{F}$ is the dimension filtration of $L\otimes_\Bbbk N $ with respect to $\mathcal{I}$.
\end{lemma}
\begin{prf}
Let $s\leq r$. We claim that 
\[
\mathscr{G}\otimes_\Bbbk \mathscr{F}: 0=L_0\otimes_{\Bbbk} N_0 \varsubsetneq \dots\varsubsetneq L_{s-1}\otimes_{\Bbbk} N_{s-1}\varsubsetneq L_s\otimes_{\Bbbk} N\varsubsetneq \dots\varsubsetneq L_{r}\otimes_{\Bbbk} N=L\otimes_{\Bbbk} N
\]
is the dimension filtration of $L\otimes_\Bbbk N$ with respect to $\mathcal{I}$. The argument for $r\leq s$ is similar.
For $i=1, \dots, s$, the strict inclusions follow from the fact that
$$\begin{array}{ll}
\cd(\mathcal{I}, L_{i-1}\otimes_{\Bbbk} N_{i-1})&=\cd(I, L_{i-1})+\cd(J, N_{i-1})\\
                                         &<\cd(I, L_{i})+\cd(J, N_{i})\\
                                         &=\cd(\mathcal{I}, L_{i}\otimes_{\Bbbk} N_{i}).\\
\end{array}$$
 Theorem  \ref{2.6}(b) yields the first and third steps, and the second step follows from our assumptions. 
For $i=s+1, \dots, r$, we have the following display
$$\begin{array}{ll}
\cd(\mathcal{I}, L_{i-1}\otimes_{\Bbbk} N)&=\cd(I, L_{i-1})+\cd(J, N)\\
                                         &<\cd(I, L_{i})+\cd(J, N)\\
                                         &=\cd(\mathcal{I}, L_{i}\otimes_{\Bbbk} N).\\
\end{array}$$
To complete our proof,  for $i=1, \dots, s$,  we need to show that $L_{i-1}\otimes_\Bbbk N_{i-1}$ is the largest submodule of $L_{i}\otimes_\Bbbk N_{i}$
for which $\cd(\mathcal{I}, L_{i-1}\otimes_\Bbbk N_{i-1})<\cd(\mathcal{I}, L_{i}\otimes_\Bbbk N_{i})$.  Moreover, for $i=s+1, \dots, r$,  a similar argument shows that  $L_{i-1}\otimes_\Bbbk N$ is the largest submodule of $L_{i}\otimes_\Bbbk N$ for which $\cd(\mathcal{I}, L_{i-1}\otimes_\Bbbk N)<\cd(\mathcal{I}, L_{i}\otimes_\Bbbk N)$.
Let $K$ be the largest submodule of $L_{i}\otimes_\Bbbk N_{i}$ such that $\cd(\mathcal{I}, K)<\cd(\mathcal{I}, L_{i}\otimes_\Bbbk N_{i})$. We want to show that $L_{i-1}\otimes_\Bbbk N_{i-1}=K$.
 By Fact \ref{4.3}, we have
  \[
 \cd(\mathcal{I}, K\oplus(L_{i-1}\otimes_\Bbbk N_{i-1}))=\max \{ \cd(\mathcal{I}, K), \cd(\mathcal{I}, L_{i-1}\otimes_\Bbbk N_{i-1})\}< \cd(\mathcal{I}, L_{i}\otimes_\Bbbk N_{i}).
  \]
  Thus, the maximality of $K$ yields $K=K\oplus (L_{i-1}\otimes_\Bbbk N_{i-1})$. Hence, $L_{i-1}\otimes_\Bbbk N_{i-1} \subseteq K \varsubsetneq L_{i}\otimes_\Bbbk N_{i}$.
 On the contrary, suppose that $L_{i-1}\otimes_\Bbbk N_{i-1} \neq  K$. Observe that
$$\begin{array}{ll}
\emptyset  \neq \Ass_{A \otimes_\Bbbk B} K/(L_{i-1}\otimes_\Bbbk N_{i-1}) & \subseteq  \Ass_{A \otimes_\Bbbk B} (L_{i}\otimes_\Bbbk N_{i})/(L_{i-1}\otimes_\Bbbk N_{i-1})\\
                                           &=  \Ass_{A \otimes_\Bbbk B} (L_i/L_{i-1})\otimes_\Bbbk (N_i/N_{i-1})\\
                                           &= \{\fp_1 \otimes_\Bbbk B+ A\otimes_\Bbbk \fp_2: \fp_1 \in \Ass_AL_i/L_{i-1} \;\; \text{and} \;\; \\
                                           &   \fp_2 \in \Ass_B  N_i/N_{i-1}\}.\\
\end{array}$$
The exact sequence
 \[
 0\to K/(L_{i-1}\otimes_\Bbbk N_{i-1})\to (L_{i}\otimes_\Bbbk N_{i})/(L_{i-1}\otimes_\Bbbk N_{i-1})
 \]
 yields the first step in this sequence. The second step is standard. Fact \ref{4.4} implies the third step.
Thus, there exists $\fq \in \Ass_{A \otimes_\Bbbk B} K/(L_{i-1}\otimes_\Bbbk N_{i-1})$ such that $\fq=\fp_1 \otimes_\Bbbk B+ A\otimes_\Bbbk \fp_2$ where $\fp_1 \in \Ass_AL_i/L_{i-1}$ and $\fp_2 \in \Ass_B  N_i/N_{i-1}$.
Observe that
$$\begin{array}{ll}
\cd(\mathcal{I}, L_{i}\otimes_\Bbbk N_{i})&=\cd(I, L_i)+\cd(J, N_i)\\
                                           &=\cd(I, A/\fp_1)+\cd(J, B/\fp_2)\\
                                           &=\cd(\mathcal{I},  A/\fp_1\otimes_\Bbbk B/\fp_2)\\
                                            &=\cd(\mathcal{I},  (A \otimes_\Bbbk B)/\fq)\\
                                           & \leq \cd(\mathcal{I}, K/(L_{i-1}\otimes_\Bbbk N_{i-1}))\\
                                           &\leq  \cd(\mathcal{I}, K).\\
\end{array}$$
The first and third steps follow from Theorem \ref{2.6}(b).  As $\cd(I, L_{i-1})<\cd(I, L_i)$, the exact sequence $0\to L_{i-1}\to L_i \to L_i/L_{i-1}\to 0$ yields 
$\cd(I, L_i)=\cd(I, L_i/L_{i-1})$ by Fact \ref{4.3}. Thus,  Fact \ref{4.2} implies the second step. The fourth step is by (\ref{2.3.1}). The exact sequence
 \[
 0\to (A \otimes_\Bbbk B)/\fq \to K/(L_{i-1}\otimes_\Bbbk N_{i-1})
 \]
yields $\cd(\mathcal{I}, (A \otimes_\Bbbk B)/\fq)\leq \cd(\mathcal{I}, K/(L_{i-1}\otimes_\Bbbk N_{i-1})$ by Fact \ref{4.3}. Thus, the fifth step follows. The exact sequence
\[
K\to K/(L_{i-1}\otimes_\Bbbk N_{i-1})\to 0
\]
yields $\cd(\mathcal{I}, K/(L_{i-1}\otimes_\Bbbk N_{i-1}))\leq \cd(\mathcal{I}, K)$ again by Fact \ref{4.3}. So, the sixth step follows.
Consequently, $\cd(\mathcal{I}, L_{i}\otimes_\Bbbk N_{i})\leq   \cd(\mathcal{I}, K)$, a contradiction. 
 Therefore, $L_{i-1}\otimes_\Bbbk N_{i-1}=K$ and so the proof is complete. 
\end{prf}
\begin{theorem}
\label{4.6}
Continue with the notation and assumptions as above. Let $\Bbbk$ be an algebraically closed field. Then the following conditions are equivalent:
 \begin{itemize}
\item[{(a)}] $L\otimes_{\Bbbk} N$ is sequentially Cohen--Macaulay with respect to $\mathcal{I}$;
\item[{(b)}] $L$ and $N$ are sequentially Cohen--Macaulay with respect to $I$ and $J$, respectively.
\end{itemize}
\end{theorem}
\begin{prf}
$(a)\Rightarrow (b)$: Assume that $L\otimes_{\Bbbk} N$ is sequentially Cohen-Macaulay with respect to $\mathcal{I}$.
 Let $\mathscr{F}: 0=L_0\varsubsetneq L_1 \varsubsetneq \dots  \varsubsetneq L_r=L$ be the dimension filtration of $L$ with respect to $I$,
 and let $\mathscr{G}: 0=N_0\varsubsetneq N_1 \varsubsetneq \dots  \varsubsetneq N_s=N$ be the dimension filtration of $N$ with respect to $J$.
 By Lemma \ref{4.5}, we may assume that $r=s$ and
\[
\mathscr{G}\otimes_\Bbbk\mathscr{F}: 0=L_0\otimes_{\Bbbk} N_0 \varsubsetneq L_1\otimes_{\Bbbk} N_1\varsubsetneq \dots\varsubsetneq L_{r}\otimes_{\Bbbk} N_{r}=L\otimes_{\Bbbk} N
\]
is the dimension filtration of $L\otimes_\Bbbk N$ with respect to $\mathcal{I}$.
 Our assumption implies that $\mathscr{G}\otimes_\Bbbk\mathscr{F}$ is the Cohen-Macaulay filtration of $L\otimes_\Bbbk N$ with respect to $\mathcal{I}$.
Observe that
$$\begin{array}{ll}
\grade(I, L_i/L_{i-1})+\grade(J, N_i/N_{i-1})&= \grade(\mathcal{I}, (L_i/L_{i-1}) \otimes_\Bbbk  (N_i/N_{i-1}))\\
                                           &=  \grade(\mathcal{I}, (L_i\otimes_\Bbbk N_i)/(L_{i-1}\otimes_\Bbbk N_{i-1}))\\
                                            &=  \cd(\mathcal{I}, (L_i\otimes_\Bbbk N_i)/(L_{i-1}\otimes_\Bbbk N_{i-1}))\\
                                            &= \cd(\mathcal{I}, (L_i/L_{i-1}) \otimes_\Bbbk (N_i/N_{i-1}))\\
                                            &=\cd(I, L_i/L_{i-1})+\cd(J, N_i/N_{i-1}).\\

\end{array}$$
The first step follows from Theorem \ref{2.6}(a). The second and fourth steps are standard. Our assumption yields the third step and the final step follows from Theorem \ref{2.6}(b).
 Consequently,  $\grade(I, L_i/L_{i-1})=\cd(I, L_i/L_{i-1})$ and $\grade(J, N_i/N_{i-1})=\cd(J, N_i/N_{i-1})$; so that $L_i/L_{i-1}$ is Cohen--Macaulay with respect to $I$ and
 $N_i/N_{i-1}$ is Cohen--Macaulay with respect to $J$ for $i=1, \dots, r$.

On the other hand, our assumptions yield
\[
\cd(I, L_i/L_{i-1})=\cd(I, L_i)<\cd(I, L_{i+1})= \cd(I, L_{i+1}/L_{i}),
\]
  and
 \[
 \cd(J, N_i/N_{i-1})=\cd(J, N_i)<\cd(J, N_{i+1})= \cd(J, N_{i+1}/N_{i}),
 \]
 for $i=1, \dots, r$.
Therefore, $L$ and $N$ are sequentially Cohen--Macaulay with respect to $I$ and $J$, respectively.

$(b)\Rightarrow (a)$:
Suppose $L$ and $N$ are sequentially Cohen--Macaulay with respect to $I$ and $J$, respectively.
Thus, there exists a Cohen-Macaulay filtration for $L$ with respect to $I$, given by $\mathscr{F}: 0=L_0\varsubsetneq L_1 \varsubsetneq \dots  \varsubsetneq L_r=L$.
Similarly, there exists a Cohen-Macaulay filtration for $N$ with respect to $J$, given by $\mathscr{G}: 0=N_0\varsubsetneq N_1 \varsubsetneq \dots  \varsubsetneq N_s=N$.
Assuming $s\leq r$, we assert that the filtration given by
\begin{eqnarray*}
0=L_0\otimes_{\Bbbk} N_1 \varsubsetneq L_1\otimes_{\Bbbk} N_1\varsubsetneq L_2\otimes_{\Bbbk} N_1\varsubsetneq
 \dots  \varsubsetneq  L_{r-s+1}\otimes_{\Bbbk} N_1  \varsubsetneq \\
  L_{r-s+2}\otimes_{\Bbbk} N_2\varsubsetneq \dots \varsubsetneq L_{r-1}\otimes_{\Bbbk} N_{s-1}\varsubsetneq L_r\otimes_{\Bbbk} N_s=L\otimes_{\Bbbk} N,
\end{eqnarray*}
is a Cohen-Macaulay filtration for $L\otimes_{\Bbbk} N$ with respect to $\mathcal{I}$. A similar argument can be made for $r\leq s$.
For $0\leq i\leq r-s$, we set
\[
M_{i,1}=(L_i\otimes_\Bbbk N_1)/(L_{i-1}\otimes_\Bbbk N_1)\cong (L_i/L_{i-1}) \otimes_\Bbbk N_1.
\]
By using Theorem \ref{2.6}, we observe that
$$\begin{array}{ll}
\grade(\mathcal{I}, M_{i, 1})&= \grade(\mathcal{I}, (L_i/L_{i-1}) \otimes_\Bbbk N_1)\\
                                                                                      &=\grade(I, L_i/L_{i-1})+\grade(J, N_1)\\
                                                                                       &=\cd(I, L_i/L_{i-1})+\cd(J, N_1)\\
                                                                                      &= \cd(\mathcal{I}, (L_i/L_{i-1}) \otimes_\Bbbk N_1)\\
                                                                                      &=\cd(\mathcal{I}, M_{i, 1}).\\
\end{array}$$
On the other hand,
$$\begin{array}{ll}
\cd(\mathcal{I}, M_{i, 1})&=\cd(I, L_i/L_{i-1})+\cd(J, N_1)\\
                                         &<\cd(I, L_{i+1}/L_i)+\cd(J, N_1)\\
                                         &=\cd(\mathcal{I}, M_{i+1,1}).\\
\end{array}$$
For $r-s+1\leq i\leq r$ and $1\leq j\leq s$, we set
\[
T_{i, j}=(L_i\otimes_\Bbbk N_j)/(L_{i-1}\otimes_\Bbbk N_{j-1})\cong (L_i/L_{i-1})\otimes_\Bbbk (N_j/N_{j-1}).
\]
Again, by using Theorem \ref{2.6}, we have the following display
$$\begin{array}{ll}
\grade(\mathcal{I}, T_{i, j})&=\grade(\mathcal{I}, (L_i/L_{i-1})\otimes_\Bbbk  (N_j/N_{j-1}))\\
                                              &=\grade(I, L_i/L_{i-1})+ \grade(J,  N_i/N_{i-1})\\
                                              &=\cd(I, L_i/L_{i-1})+\cd(J,  N_j/N_{j-1})\\
                                              &=\cd(\mathcal{I}, (L_i/L_{i-1})\otimes_\Bbbk (N_j/N_{j-1}))\\.
                                              &=\cd(\mathcal{I}, T_{i, j}).\\

\end{array}$$
On the other hand,
$$\begin{array}{ll}
\cd(\mathcal{I}, T_{i, j})&=\cd(I, L_i/L_{i-1})+\cd(J, N_j/N_{j-1})\\
                                         &<\cd(I, L_{i+1}/L_i)+\cd(J, N_{j+1}/N_{j})\\
                                         &=\cd(\mathcal{I}, T_{i+1, j+1}).\\
\end{array}$$
Therefore, the proof is complete.
\end{prf}
As a consequence of Theorem \ref{4.6}, we have the following result:
\begin{corollary}
\label{4.7}
Continue with the notation and assumptions as above. Let $\Bbbk$ be an algebraically closed field. Then the following conditions are equivalent:
 \begin{itemize}
\item[{(a)}] $L\otimes_{\Bbbk} N$ is sequentially Cohen--Macaulay with respect to $I \otimes_\Bbbk B$ (or $A \otimes_\Bbbk J$).
\item[{(b)}] $L$ is sequentially Cohen--Macaulay with respect to $I$ (or $N$ is sequentially Cohen--Macaulay with respect to $J$).
\end{itemize}
\end{corollary}
As a further implication of Theorem \ref{4.6}, we obtain the following result. A graded variant of this fact can be found in \cite[Theorem 2.11]{STY}.
\begin{corollary}
\label{4.8}
Let $\Bbbk$ be an algebraically closed field. Let $(A, \fm)$ and $(B, \fn)$ be two Noetherian local algebras over a field $\Bbbk$ such that $A \otimes_\Bbbk B$ is Noetherian.  We also assume 
that $\Bbbk$ is the common residue field of both $A$ and $B$, and either $A$ or $B$ is algebraic over $\Bbbk$.  
  Then the following conditions are equivalent:
 \begin{itemize}
\item[{(a)}] $L\otimes_{\Bbbk} N$ is sequentially Cohen--Macaulay $A\otimes_\Bbbk B$-module;
\item[{(b)}] $L$ and $N$ are sequentially Cohen--Macaulay over $A$ and $B$, respectively.
\end{itemize}
\end{corollary}
We recall the following fact:
\begin{fact}
\label{4.9} {\em
Let $R$ be a Noetherian ring, $\fa$ an ideal of $R$, and $M$ a finitely generated $R$-module. Assume that $M$ is sequentially Cohen--Macaulay with respect to $\fa$ with the Cohen--Macaulay filtration $\mathcal{F}$: $0=M_0\varsubsetneq M_1 \varsubsetneq  \dots \varsubsetneq M_r=M$ with respect to $\fa$. Then
\begin{itemize}
\item[{(a)}] $ \text{H}^{c_i}_\fa(M)\cong \text{H}^{c_i}_\fa(M_i)\cong \text{H}^{c_i}_\fa(M_i/M_{i-1})$ where $c_i=\cd(\fa, M_i)$ for $i=1, \dots, r$ and $\text{H}^{k}_\fa(M)=0$ for $k\not \in \{ c_1, \dots, c_r\}$, see \cite[Lemma 3.2]{R3}.
\item[{(b)}]  $\grade(\fa, M_i)=\grade(\fa, M)$ for $i=1, \dots, r$, see \cite[Fact 3.3]{R3}.
\end{itemize}
We remark that the proof of these facts is based on the same method as in the proofs of \cite[Lemma 3.2 and Fact 3.3]{R3} by replacing the ring $R$ and the general ideal $\fa$ with the polynomial ring $S$ and the ideal $Q$, respectively.
}
\end{fact}

\begin{proposition}
\label{4.10}
Let $L$ and $N$ are sequentially Cohen--Macaulay with respect to $I$ and $J$, respectively.
Assume that $\grade(I, L)>0$ and $\grade(J, N)>0$.  Then
\[
\f_{\mathcal{I}}(L\otimes_{\Bbbk} N)=\f_I(L)+\f_J(N).
\]
\end{proposition}
\begin{prf}
Based on our assumptions, there is a Cohen-Macaulay filtration for $L$ with respect to $I$, represented as $\mathscr{F}: 0=L_0\varsubsetneq L_1 \varsubsetneq \dots  \varsubsetneq L_r=L$, and a Cohen-Macaulay filtration for $N$ with respect to $J$, represented as $\mathscr{G}: 0=N_0\varsubsetneq N_1 \varsubsetneq \dots  \varsubsetneq N_s=N$. 
It can be observed that
\[
\begin{aligned}
\f_I(L) &= \cd(I, L_1) \\
&= \grade(I, L_1) \\
&= \grade(I, L).
\end{aligned}
\]
The first step in this sequence follows from Fact \ref{4.9}(a). Since $L_1$ is Cohen--Macaulay with respect to $I$, the second step follows. The third step follows from Fact \ref{4.9}(b). Similarly, $\f_J(N)=\grade(J, N)$. 
The assertion now follows from Theorem \ref{3.2}.
\end{prf}
\begin{corollary}
\label{4.11}
Let $(A, \fm)$ and $(B, \fn)$ be two Noetherian local algebras over a field $\Bbbk$ such that $A \otimes_\Bbbk B$ is Noetherian.  We also assume 
that $\Bbbk$ is the common residue field of both $A$ and $B$, and either $A$ or $B$ is algebraic over $\Bbbk$.   Let $L$ and $N$ be sequentially Cohen--Macaulay modules with $\depth_A L >0$ and $\depth_B N>0$.  Then
\[
\f_{\mathfrak{M}}(L\otimes_{\Bbbk} N)=\f_\fm(L)+\f_\fn(N).
\]
\end{corollary}
\begin{remark}
\label{4.12}{\em
Example \ref{3.7} illustrates that the conditions $\depth_A L >0$ and $\depth_B N>0$ in Corollary \ref{4.11} are necessary. In fact, since any $R$-module $M$ of dimension less than or equal to one is sequentially Cohen--Macaulay, we have that $L$ and $N$ are sequentially Cohen--Macaulay. However, it should be noted that $\depth_A L=0$, so the equality $\f_{\mathfrak{M}}(L\otimes_{\Bbbk} N)=\f_\fm(L)+\f_\fn(N)$ does not hold.
}
\end{remark}

\begin{center}
{Acknowledgment}
\end{center}
\hspace*{\parindent} The author would like to thank the referee for their careful review and helpful comments.

\end{document}